\documentclass[12pt,smallextended]{svjour3}
\smartqed
\usepackage{lineno}
\modulolinenumbers[5]

\journalname{Journal of XXX}
\usepackage{amssymb,amsmath,amsfonts}
\usepackage[bookmarksnumbered, colorlinks, urlcolor=blue, linkcolor=blue, citecolor=red, plainpages=false, pdfstartview=FitW]{hyperref}
\usepackage{array}
\usepackage{bm}
\usepackage{epsfig}
\usepackage{mathptmx}
\usepackage{diagbox}
\usepackage{graphicx}
\usepackage{color}
\usepackage{caption}
\usepackage{wrapfig}
\usepackage{graphics}
\usepackage{enumerate}
\usepackage{amsxtra}
\usepackage{latexsym}
\usepackage{multirow}
\usepackage{slashbox}
\usepackage{subfigure}
\usepackage{epstopdf}
\usepackage{pdfpages}
\usepackage{algorithm}
\usepackage{algorithmic}
\usepackage{url}
\newcommand{\thmlist}{
\begin{list}{Step 1}
{\setlength{\leftmargin}{0.6 in}\setlength{\labelwidth} {0.5 in}}}

\newcommand{\alglist}{
\begin{list}{Step 1}
{\setlength{\leftmargin}{1.1 in} \setlength{\labelwidth}{1.0 in}}}

 \renewcommand{\proof} {\noindent {\bf Proof.} \quad}

\renewcommand{\subtitle}[1]{\color{blue}}

\begin{document}

%\begin{frontmatter}

\title{Explicit continuation methods with L-BFGS updating formulas
for linearly constrained optimization problems}
\titlerunning{Explicit Continuation methods with L-BFGS updating formulas}
\author{Xin-long Luo\textsuperscript{$\ast$} \and Jia-hui Lv \and Hang Xiao}
\authorrunning{Luo, Lv and Xiao}
%\authorrunning{Short form of author list} % if too long for running head

\institute{Xin-long Luo
             \at
              Corresponding author. School of Artificial Intelligence, \\
Beijing University of Posts and Telecommunications, P. O. Box 101, \\
Xitucheng Road  No. 10, Haidian District, 100876, Beijing China\\
             \email{luoxinlong@bupt.edu.cn}            %  \\
        \and
        Jia-hui Lv \at
        School of Artificial Intelligence, \\
Beijing University of Posts and Telecommunications, P. O. Box 101, \\
Xitucheng Road  No. 10, Haidian District, 100876, Beijing China\\
             \email{jhlv@bupt.edu.cn}
           \and
           Hang Xiao
     \at
     School of Artificial Intelligence, \\
     Beijing University of Posts and Telecommunications, P. O. Box 101, \\
     Xitucheng Road  No. 10, Haidian District, 100876, Beijing China \\
     \email{xiaohang0210@bupt.edu.cn}
}

\date{Received: date / Accepted: date}
% The correct dates will be entered by the editor
\maketitle

\begin{abstract}
  This paper considers an explicit continuation method with the trusty 
  time-stepping scheme and the limited-memory BFGS (L-BFGS) 
  updating formula (Eptctr) for the linearly constrained optimization 
  problem. At every iteration, Eptctr only involves three pairs of the 
  inner product of vector and one
  matrix-vector product, other than the traditional and representative
  optimization method such as the sequential quadratic programming (SQP)
  or the latest continuation method such as Ptctr \cite{LLS2020}, which
  needs to solve a quadratic programming subproblem (SQP) or a linear
  system of equations (Ptctr). Thus, Eptctr can save much more computational
  time than SQP or Ptctr. Numerical results also show that the consumed time
  of EPtctr is about one tenth of that of Ptctr or one fifteenth to 0.4
  percent of that of SQP. Furthermore, Eptctr can save the storage space of an
  $(n+m) \times (n+m)$ large-scale matrix, in comparison to SQP. The required
  memory of Eptctr is about one fifth of that of SQP. Finally, we also give the
  global convergence analysis of the new method under the standard assumptions.
\end{abstract}

% keywords here, in the form: keyword \sep keyword


\vskip 2mm

%\textbf{AMS subject classifications.} 65H17 \and 65J15 \and 65K05 \and 65L05

%\end{frontmatter}

% \linenumbers
% main text

\section{Introduction} \label{SUBINT}

\vskip 2mm

In this article, we consider the following linearly  equality-constrained
optimization problem
\begin{align}
  &\min_{x \in \Re^n} \; f(x)  \nonumber \\
  &\text{subject to} \; \; Ax = b,   \label{LEQOPT}
\end{align}
where matrix $A \in \Re^{m \times n}$ and vector $b \in \Re^{m}$
may have random noise. This problem has many applications in engineering
fields such as the visual-inertial navigation of an unmanned aerial
vehicle maintaining the horizontal flight \cite{CMFO2009,LLS2020},
and there are many practical methods to solve it such as the sequential
quadratic programming (SQP) method \cite{Bertsekas2018,NW1999} or the penalty
function method \cite{FM1990}.

\vskip 2mm

For the constrained optimization problem \eqref{LEQOPT}, the continuation
method \cite{AG2003,CKK2003,Goh2011,KLQCRW2008,Pan1992,Tanabe1980} is another
method other than the traditional optimization method such as SQP or the penalty
function method. The advantage of the continuation method over the SQP method
is that the continuation method is capable of finding many local optimal points
of the non-convex optimization problem by tracking its trajectory, and it is
even possible to find the global optimal solution
\cite{BB1989,Schropp2000,Yamashita1980}. However, the computational
efficiency of the continuation method  may be higher than that of SQP. Recently,
Luo, Lv and Sun \cite{LLS2020} give a continuation method with the trusty time-stepping
scheme and its consumed time is about one fifth of that of SQP for the linearly
constrained optimization problem \eqref{LEQOPT}. Their method only needs to
solve a linear system of equations with an $n \times n$ symmetric definite
coefficient matrix at every iteration, which involves about $\frac{1}{3}n^{3}$
flops. SQP needs to solve a linear system of equations with an
$(m+n) \times (m+n)$ coefficient matrix, which involves about $\frac{2}{3}(m+n)^{3}$
flops. In order to improve the computational efficiency further and save the storage
of the continuation method \cite{LLS2020} for the large-scale optimization problem,
we consider a special limited-memory BFGS updating formula and the trusty
time-stepping scheme in this article.

\vskip 2mm

The rest of the paper is organized as follows. In section 2, we give a new
continuation method with the trusty time-stepping scheme and the L-BFGS updating 
formula for the linearly equality-constrained optimization problem
\eqref{LEQOPT}. In section 3, we analyze the global convergence of this new
method. In section 4, we report some promising numerical results of the new
method, in comparison to the traditional optimization method (SQP) and the
latest continuation method (Ptctr) for some large-scale problems. Finally, we
give some discussions and conclusions in section 5.

\vskip 2mm

\section{The explicit continuation method with L-BFGS updating formulas}

\vskip 2mm

In this section, we construct an explicit continuation method with the adaptive
time-stepping scheme based on the trust-region updating strategy \cite{Yuan2015}
for the linearly equality-constrained optimization problem \eqref{LEQOPT}.
Firstly, we construct a generalized projected gradient flow based on the KKT
conditions of linearly constrained optimization problem. Then, in order to
efficiently follow the generalized gradient flow, we construct an explicit
continuation method with an adaptive time-stepping scheme for this special
ordinary differential equations (ODEs). Furthermore, we give a preprocessing
method for the infeasible initial point.

\vskip 2mm

\subsection{The generalized projected gradient flow}

\vskip 2mm

For the linearly constrained optimization problem \eqref{LEQOPT}, it is well known
that its optimal solution $x^{\ast}$ needs to satisfy the Karush-Kuhn-Tucker
conditions (p. 328, \cite{NW1999}) as follows:
\begin{align}
  \nabla_{x} L(x, \, \lambda) &= \nabla f(x) + A^{T} \lambda = 0,
    \label{FOKKTG} \\
  Ax - b & = 0,             \label{FOKKTC}
\end{align}
where the Lagrangian function $L(x, \, \lambda)$ is defined by
\begin{align}
  L(x, \, \lambda) = f(x) + \lambda^{T}(Ax-b).
      \label{LAGFUN}
\end{align}
Similarly to the method of the negative gradient flow for the unconstrained
optimization problem \cite{LKLT2009}, from the first-order necessary conditions
\eqref{FOKKTG}-\eqref{FOKKTC}, we can construct a dynamical system of
differential-algebraic equations for problem \eqref{LEQOPT}
\cite{LL2010,Luo2012,LLW2013,Schropp2003} as follows:
\begin{align}
    & \frac{dx}{dt} = - \nabla L_{x}(x, \, \lambda)
      = -\left(\nabla f(x) + A^{T} \lambda \right),  \label{DAGF} \\
    & Ax - b = 0.                      \label{LACON}
\end{align}

\vskip 2mm

By differentiating the algebraic constraint \eqref{LACON} with respect to $t$
and replacing it into the differential equation \eqref{DAGF}, we obtain
\begin{align}
  A\frac{dx}{dt} = - A \left(\nabla f(x) + A^{T} \lambda \right)
  = - A \nabla f(x) - AA^{T} \lambda = 0.    \label{DIFALGC}
\end{align}
If we assume that matrix $A$ has full row rank further, from equation
\eqref{DIFALGC}, we obtain
\begin{align}
   \lambda = - \left(AA^{T} \right)^{-1} A \nabla f(x). \label{LAMBDA}
\end{align}
By replacing $\lambda$ of equation \eqref{LAMBDA} into equation \eqref{DAGF},
we obtain the dynamical system of ordinary differential equations (ODEs) as
follows:
\begin{align}
  \frac{dx}{dt} = - \left( I - A^{T} \left(AA^{T}\right)^{-1}A\right)
  \nabla f(x). \label{ODGF}
\end{align}
Thus, we also obtain the projected gradient flow for the constrained optimization
problem \cite{Tanabe1980}.

\vskip 2mm

For convenience, we denote the projection matrix $P$ as
\begin{align}
  P  = I - A^{T} \left(AA^{T}\right)^{-1}A.  \label{PROMAT}
\end{align}
It is not difficult to verify $P^{2} = P$ and $AP = 0$. That is to say, $P$
is a symmetric projection matrix and its eigenvalues are 0 or 1. From Theorem 2.3.1
in p. 73 of \cite{GV2013}, we know that its matrix 2-norm is
\begin{align}
     \|P\| = 1. \label{MATNP}
\end{align}
We denote $P^{+}$ as the Moore-Penrose generalized inverse of $P$ (p. 11, \cite{SY2006}).
Since $P$ is symmetric and its eigenvalues are 0 or 1, it is not difficult to verify
\begin{align}
      P^{+} = P.     \label{GINVP}
\end{align}
Thus, for a full rank matrix $B \in \Re^{n \times n}$, we obtain the generalized inverse
$(PB)^{+}$ of $PB$ as follows:
\begin{align}
      (PB)^{+} = B^{+}P^{+} = B^{-1}P. \label{GINVBP}
\end{align}

\vskip 2mm

Similarly to the generalized gradient flow for an unconstrained optimization problem
(p. 361, \cite{HM1996}), from the projected gradient flow \eqref{ODGF}, we can
construct the generalized projected gradient flow for the constrained optimization
problem \eqref{LEQOPT} as follows:
\begin{align}
    \frac{dx}{dt} = - (PH(x)P)P\nabla f(x) ,
    \; x(0) = x_{0}, \label{GPGF}
\end{align}
where $H(x)$ is a symmetric positive definite matrix for any $x \in \Re^{n}$.
Here, $H(x)$ may be selected as the inverse of the Hessian matrix
$\nabla^{2} f(x)$ of $f(x)$ and $PH(x)P$ can be regarded as a pre-conditioner
of $P\nabla f(x)$ to mitigate the stiffness of the ODEs \eqref{GPGF}.
Consequently, we can adopt the explicit numerical method to compute the
trajectory of the ODEs \eqref{GPGF} efficiently \cite{LXL2020,LY2021}.

\begin{remark}
If $x(t)$ is the solution of the ODEs \eqref{GPGF}, it is not difficult to verify
that $x(t)$ satisfies $A (dx/dt) = 0$. That is to say, if the initial point
$x_{0}$ satisfies $Ax_{0} = b$, the solution $x(t)$ of the generalized projected
gradient flow \eqref{GPGF} also satisfies $Ax(t) = b, \; \forall t \ge 0$. This
property is very useful when we construct a structure-preserving algorithm
\cite{HLW2006,Simos2013} to follow the trajectory of the ODEs \eqref{GPGF} and obtain
its equilibrium point $x^{\ast}$.
\end{remark}

\vskip 2mm

\begin{remark}
If we assume that $x(t)$ is the solution of the ODEs \eqref{GPGF}, from equations
\eqref{PROMAT}-\eqref{MATNP} and the positive definite property of $H(x)$,
we obtain
\begin{align}
 \frac{df(x)}{dt} & = \left(\nabla f(x)\right)^{T} \frac{dx}{dt}
 = - (\nabla f(x))^{T} PH(x)P \nabla f(x) \nonumber \\
 & = - (P \nabla f(x))^{T}H(x)(P \nabla f(x)) \le 0.
 \nonumber
\end{align}
That is to say, $f(x)$ is monotonically decreasing along the solution curve $x(t)$
of the dynamical system \eqref{GPGF}. Furthermore, the solution $x(t)$ converges
to $x^{\ast}$ when $f(x)$ is lower bounded and $t$ tends to infinity
\cite{HM1996,Schropp2000,Tanabe1980}, where $x^{\ast}$ satisfies the first-order
Karush-Kuhn-Tucker conditions \eqref{FOKKTG}-\eqref{FOKKTC}. Thus, we can follow
the trajectory $x(t)$ of the ODEs \eqref{GPGF} to obtain its equilibrium point $x^{\ast}$,
which is also one saddle point of the original optimization problem \eqref{LEQOPT}.
\end{remark}

\vskip 2mm

\subsection{The explicit continuation method} \label{SUBSICM}

\vskip 2mm

The solution curve of the degenerate ordinary differential equations is not
efficiently followed on an infinite interval by the traditional ODE method
\cite{AP1998,BCP1996,LF2000}, so one needs to construct the
particular method for this problem \eqref{GPGF}. We apply the first-order implicit
Euler method \cite{SGT2003} to the ODEs \eqref{GPGF}, then we obtain
\begin{align}
      x_{k+1} =  x_{k} - \Delta t_{k} (PH(x_{k}))(P\nabla f(x_{k})),
      \label{IEGPGF}
\end{align}
where $\Delta t_k$ is the time-stepping size.

\vskip 2mm

Since the system of equations \eqref{IEGPGF} is a nonlinear system which is not
directly solved, we seek for its explicit approximation formula.
We denote $s_{k} = x_{k+1} - x_{k}$. By using the first-order Taylor expansion,
we have the linear approximation $\nabla f(x_{k}) + \nabla^{2} f(x_{k})s_{k}$ of
$\nabla f(x_{k+1})$. By substituting it into equation \eqref{IEGPGF} and using the
zero-order approximation $H(x_{k})$ of $H(x_{k+1})$, we have
\begin{align}
      s_{k} &\approx - \Delta t_{k} (PH(x_{k}))\left(P\left(\nabla f(x_{k}) +
     \nabla^{2} f(x_{k})s_{k}\right)\right) \nonumber \\
     & = - \Delta t_{k} (PH(x_{k}))(P \nabla f(x_{k}))
     - \Delta t_{k} P(H(x_{k})P)(P\nabla^{2}f(x_{k}))s_{k}. \label{AIEGPGF}
\end{align}
From equation \eqref{IEGPGF} and $P^{2} = P$, we have $Ps_{k} = s_{k}$.
Let $H(x_{k}) = (\nabla^{2} f(x_{k}))^{-1}$. Then, we have $H(x_{k})P =
(P\nabla^{2}f(x_{k}))^{+}$. Thus, we regard
\begin{align}
      P (H(x_{k})P)(P\nabla^{2}f(x_{k})Ps_{k})
      = P (P\nabla^{2}f(x_{k}))^{+} (P\nabla^{2}f(x_{k})) Ps_{k} \approx Ps_{k}
      = s_{k}. \label{AGINVPHP}
\end{align}
By substituting it into equation \eqref{AIEGPGF}, we obtain the explicit
continuation method as follows:
\begin{align}
      & s_{k} = - \frac{\Delta t_{k}}{1 + \Delta t_{k}} (PH_{k})(Pg_{k}),
      \label{EXPTC} \\
      & x_{k+1} = x_{k} + s_{k}, \label{XK1}
\end{align}
where $g_{k} = \nabla f(x_{k})$ and $H_{k} = (\nabla^{2} f(x_{k}))^{-1}$ or
its quasi-Newton approximation in the projective space $S_{p}^{k} = \{x: \;
x = x_{k} + Pd, \, d \in \Re^{n}\}$.

\vskip 2mm

If we let the projection matrix $P = I$, the formula \eqref{EXPTC} is equivalent
to the explicit continuation method given by Luo, Xiao and Lv \cite{LXL2020}
for nonlinear equations. The explicit continuation method \eqref{EXPTC}-\eqref{XK1}
is similar to the projected damped Newton method if we let $\alpha_{k} =
\Delta t_k/(1+\Delta t_k)$ in equation \eqref{EXPTC}. However,
from the view of the ODE method, they are different. The projected damped
Newton method is obtained by the explicit Euler scheme applied
to the generalized projected gradient flow \eqref{GPGF}, and its time-stepping size
$\alpha_k$ is restricted by the numerical stability \cite{SGT2003}.
That is to say, the large time-stepping size $\alpha_{k}$ can not be adopted in the
steady-state phase.

\vskip 2mm

The explicit continuation method \eqref{EXPTC}-\eqref{XK1} is obtained by the 
implicit Euler approximation method applied to the generalized projected
gradient flow \eqref{GPGF}, and its time-stepping size $\Delta t_k$ is not restricted
by the numerical stability. Therefore, the large time-stepping size can be adopted
in the steady-state phase for the explicit continuation method
\eqref{EXPTC}-\eqref{XK1}, and it mimics the Newton method near the equilibrium
solution $x^{\ast}$ such that it has the fast local convergence rate. The most
of all, the new step size $\alpha_{k} = \Delta t_{k}/(\Delta t_{k} + 1)$ is
favourable to adopt the trust-region updating technique for adaptively adjusting
the time-stepping size $\Delta t_{k}$ such that the explicit continuation method
\eqref{EXPTC}-\eqref{XK1} accurately tracks the trajectory of the generalized
projected gradient flow in the transient-state phase and achieves the fast
convergence rate near the equilibrium point $x^{\ast}$.

\vskip 2mm

\begin{remark} \label{REMPLS}
From equation \eqref{EXPTC} and the property $AP = 0$ of the
projected matrix $P$, it is not difficult to verify $As_{k} = 0$. Thus, if
the initial point $x_{0}$ satisfies the linear constraint $Ax_{0} = b$, the
point $x_{k}$ also satisfies the linear constraint $Ax_{k} = b$. That is to
say, the explicit continuation method \eqref{EXPTC}-\eqref{XK1} is a
structure-preserving method.
\end{remark}

\vskip 2mm

\subsection{The L-BFGS quasi-Newton updating formula}

\vskip 2mm

For the large-scale problem, the numerical evaluation of the Hessian matrix
$\nabla^{2}f(x_{k})$ consumes much time and stores an $n \times n$ matrix.
In order to overcome these two shortcomings, we use the L-BFGS
quasi-Newton formula (\cite{BNY1987,Goldfarb1970} or pp. 222-230, \cite{NW1999})
to approximate the generalized inverse $H(x_{k})P$ of $P\nabla^{2}f(x_{k})$.
Recently, Ullah, Sabi and Shah \cite{USS2020} give an efficient L-BFGS updating
formula for the system of monotone nonlinear equations. Here, in order to suit
the generalized projected gradient flow \eqref{GPGF}, we revise their L-BFGS
updating formula as
\begin{align}
      H_{k+1} = \begin{cases}
                      I  - \frac{y_{k}s_{k}^{T} + s_{k}y_{k}^{T}}{y_{k}^{T}s_{k}}
                      + 2\frac{y_{k}^{T}y_{k}}{(y_{k}^{T}s_{k})^{2}} s_{k}s_{k}^{T},
                      \; \text{if} \; |s_{k}^{T}y_{k}| > \theta \|s_{k}\|^{2},  \\
                      I, \; \text{otherwise}.
                \end{cases}
      \label{LBFGS}
\end{align}
where $s_{k} = x_{k+1} - x_{k}, \; y_{k} = P\nabla f(x_{k+1}) - P\nabla f(x_{k})$
and $\theta$ is a small positive constant such as $\theta = 10^{-6}$.
The initial matrix $H_{0}$ can be simply selected by the identity matrix. When
$|s_{k}^{T}y_{k}| \ge \theta \|s_{k}\|^{2}$, from equation \eqref{LBFGS}, it is
not difficult to verify
\begin{align}
       H_{k+1}y_{k} = \frac{y_{k}^{T}y_{k}}{y_{k}^{T}s_{k}} s_{k}. \nonumber
\end{align}
That is to say, $H_{k+1}$ satisfies the scaling quasi-Newton property. By using the
Sherman-Morrison-Woodburg formula, from equation \eqref{LBFGS}, when
$|s_{k}^{T}y_{k}| \ge \theta \|s_{k}\|^{2}$, we have
\begin{align}
      B_{k+1} = H_{k+1}^{-1} = I - \frac{s_{k}s_{k}^{T}}{s_{k}^{T}s_{k}}
        + \frac{y_{k}y_{k}^{T}}{y_{k}^{T}y_{k}}.\nonumber
\end{align}

\vskip 2mm

The L-BFGS updating formula \eqref{LBFGS} has some nice properties such as the
symmetric positive definite property and the positive lower bound of its eigenvalues.

\vskip 2mm

\begin{lemma} \label{LEMHLB}
Matrix $H_{k+1}$ defined by equation \eqref{LBFGS} is symmetric positive definite
and its eigenvalues are greater than $1/2$.
\end{lemma}

\vskip 2mm

\proof (i) For any nonzero vector $z \in \Re^{n}$, from equation \eqref{LBFGS},
when $|s_{k}^{T}y_{k}| > \theta \|s_{k}\|^{2}$, we have
\begin{align}
       & z^{T}H_{k+1}z = \|z\|^{2} - 2{(z^{T}y_{k})(z^{T}s_{k})}/{y_{k}^{T}s_{k}}
       + 2(z^{T}s_{k})^{2}{\|y_{k}\|^{2}}/{(y_{k}^{T}s_{k})^{2}} \nonumber \\
       & \quad =  \left(\|z\| -
         \left|{z^{T}s_{k}}/{y_{k}^{T}s_{k}}\right|\|y_{k}\|\right)^{2}
          + 2 \|z\| \left|{z^{T}s_{k}}/{y_{k}^{T}s_{k}}\right|\|y_{k}\|
          \nonumber \\
       & \quad \quad  - 2 {(z^{T}y_{k})(z^{T}s_{k})}/{y_{k}^{T}s_{k}}
          +\|y_{k}\|^{2} (z^{T}s_{k}/y_{k}^{T}s_{k})^{2}
          \ge 0.  \label{ZTHK1Z}
\end{align}
In the last inequality of equation \eqref{ZTHK1Z}, we use the Cauchy-Schwartz
inequality $\|z^{T}y\| \le \|z\|\|y_{k}\|$ and its equality holds if only if
$z = t y_{k}$. When $z = t y_{k}$, from equation \eqref{ZTHK1Z}, we have
$z^{T}H_{k+1}z = t^{2}\|y_{k}\|^{2} = \|z\|^{2}> 0$. When $z^{T}s_{k} = 0$,
from equation \eqref{ZTHK1Z}, we also have $z^{T}H_{k+1}z = \|z\|^{2} > 0$.
Therefore, we conclude that $H_{k+1}$ is a symmetric positive definite matrix
when $|s_{k}^{T}y_{k}| > \theta \|s_{k}\|^{2}$. From equation \eqref{LBFGS},
We apparently conclude that $H_{k+1}$ is a symmetric positive definite matrix
since $H_{k+1} = I$ when $|s_{k}^{T}y_{k}| \le \theta \|s_{k}\|^{2}$.

\vskip 2mm

(ii) It is not difficult to know that it exists at least $n-2$ linearly independent
vectors $z_{1}, \, z_{2}, \, \ldots, \, z_{n-2}$ such that $z_{i}^{T}s_{k} = 0, \,
z_{i}^{T}y_{k} = 0 \, ( i = 1 :(n-2))$ hold. That is to say, matrix $H_{k+1}$ defined
by equation \eqref{LBFGS} has at least $(n-2)$ linearly independent eigenvectors
whose corresponding eigenvalues are 1. We denote the other two eigenvalues of
$H_{k+1}$ as $\mu_{i}^{k+1} \, (i = 1:2)$ and their corresponding eigenvalues as
$p_{1}$ and $p_{2}$, respectively. Then, from equation \eqref{LBFGS}, we know that
the eigenvectors $p_{i} \, (i = 1:2)$ can be represented as $p_{i} = y_{k}
+ \beta_{i} s_{k}$ when $\mu_{i}^{k+1} \neq 1 \, (i = 1:2)$. From equation
\eqref{LBFGS} and $H_{k+1}p_{i} = \mu_{i}^{k+1} p_{i} \, (i = 1:2)$, we have
\begin{align}
     - \left(\mu_{i}^{k+1} + \beta_{i}\frac{s_{k}^{T}s_{k}}{s_{k}^{T}y_{k}}\right)y_{k}
     + \left(\frac{y_{k}^{T}y_{k}}{y_{k}^{T}s_{k}}
     + 2 \beta_{i} \frac{(y_{k}^{T}y_{k})(s_{k}^{T}s_{k})}{(y_{k}^{T}s_{k})^{2}}
     - \mu_{i}^{k+1} \beta_{i}\right)s_{k} = 0, \; i = 1:n. \label{EIGVAS}
\end{align}

\vskip 2mm

When $y_{k} = ts_{k}$, from equation \eqref{LBFGS}, we have $H_{k+1} = I$. In this
case, we conclude that the eigenvalues of $H_{k+1}$ are greater than $1/2$. When
vectors $y_{k}$ and $s_{k}$ are linearly independent, from equation \eqref{EIGVAS},
we have
\begin{align}
      & \mu_{i}^{k+1} + \beta_{i} {s_{k}^{T}s_{k}}/{s_{k}^{T}y_{k}} = 0,
      \nonumber \\
      &
      {y_{k}^{T}y_{k}}/{y_{k}^{T}s_{k}}
      + 2 \beta_{i}{(y_{k}^{T}y_{k})(s_{k}^{T}s_{k})}/{(y_{k}^{T}s_{k})^{2}}
     - \mu_{i}^{k+1}\beta_{i} = 0, \; i = 1:n. \nonumber
\end{align}
That is to say, $\mu_{i}^{k+1} \, (i = 1:2)$ are the two solutions of the following
equation:
\begin{align}
      \mu^{2} - 2\mu (y_{k}^{T}y_{k})(s_{k}^{T}s_{k})/(s_{k}^{T}y_{k})^{2}
      + (y_{k}^{T}y_{k})(s_{k}^{T}s_{k})/(s_{k}^{T}y_{k})^{2} = 0. \label{QUDEDQ}
\end{align}
Consequently, from equation \eqref{QUDEDQ}, we obtain
\begin{align}
      \mu_{1}^{k+1} + \mu_{2}^{k+1} = 2 (y_{k}^{T}y_{k})(s_{k}^{T}s_{k})/(s_{k}^{T}y_{k})^{2},
      \;
      \mu_{1}^{k+1}\mu_{2}^{k+1} = (y_{k}^{T}y_{k})(s_{k}^{T}s_{k})/(s_{k}^{T}y_{k})^{2}.
      \label{ROOTQEQ}
\end{align}
From equation \eqref{ROOTQEQ}, it is not difficult to obtain
\begin{align}
       {1}/{\mu_{1}^{k+1}} + {1}/{\mu_{2}^{k+1}} = 2, \;  \mu_{i}^{k+1} > 0,
       \; i = 1:2.    \label{REROOT}
\end{align}
Therefore, from equation \eqref{REROOT}, we conclude that
$\mu_{i}^{k+1} > \frac{1}{2}  \, (i = 1:2)$. Consequently, the eigenvalues of $H_{k+1}$
are greater than 1/2. \qed

\vskip 2mm

If $s_{k-1}$ is obtained from the explicit continuation method \eqref{EXPTC},
we have $Ps_{k-1} = s_{k-1}$ since $P^{2} = P$. By combining it with the L-BFGS
updating formula \eqref{LBFGS}, the explicit continuation method
\eqref{EXPTC}-\eqref{XK1} can be simplified by
\begin{align}
      & d_{k} = \begin{cases}
                   - p_{g_{k}}, \; \text{if}
                   \; |s_{k-1}^{T}y_{k-1}| \le \theta \|s_{k-1}\|^{2}, \\
                   - p_{g_{k}} + \frac{y_{k-1}(s_{k-1}^{T}p_{g_{k}})
                   + s_{k-1}(y_{k-1}^{T}p_{g_{k}})}{y_{k-1}^{T}s_{k-1}}
                   - 2 \frac{\|y_{k-1}\|^{2}(s_{k-1}^{T}p_{g_{k}})}
                  {(y_{k-1}^{T}s_{k-1})^{2}}s_{k-1}, \, \text{otherwise},
                \end{cases}         \label{SEXPTC} \\
      & s_{k} = \frac{\Delta t_{k}}{1 + \Delta t_{k}}d_{k}, \;
       x_{k+1} = x_{k} + s_{k}, \label{SXK1}
\end{align}
where $p_{g_{k}} = Pg_{k} = P \nabla f(x_{k})$ and $y_{k-1}
= P\nabla f(x_{k}) - P \nabla f(x_{k-1})$. Thus, it does not need to store
the matrix $H_{k}$ in practical computation. Furthermore, it only requires
three pairs of the inner product of vector and one matrix-vector product
($p_{g_{k}} = Pg_{k}$) to obtain the trial step $s_{k}$ and involves in
$O((n-m)n)$ flops when we use the QR decomposition or the
singular value decomposition to obtain the projection matrix $P$ in
subsection \ref{SUBSECDEG}.

\vskip 2mm

\subsection{The treatments of infeasible initial points and projection matrices}
\label{SUBSECDEG}

\vskip 2mm

We need to compute the projected gradient $Pg_{k}$ at every iteration in the updating
formula \eqref{SEXPTC}. In order to reduce the computational complexity, we use
the QR factorization (pp.276-278, \cite{GV2013}) to factor $A^{T}$ into a product
of an orthogonal matrix $Q \in \Re^{n \times n}$ and an upper triangular matrix
$R \in \Re^{n \times m}$:
\begin{align}
      A^{T} = QR = \begin{bmatrix} Q_{1} | Q_{2} \end{bmatrix}
      \begin{bmatrix} R_{1} \\ 0     \end{bmatrix}, \label{ATQR}
\end{align}
where $Q_{1} = Q(1:n, \, 1:m), \; Q_{2} = Q(1:n, \, m+1:n)$,
$R_{1} = R(1:r, \, 1:m)$ is upper triangular and nonsingular. Then, from 
equations \eqref{PROMAT}, \eqref{ATQR}, we simplify the projection matrix $P$ as
\begin{align}
      P = I - Q_{1}Q_{1}^{T} = Q_{2}Q_{2}^{T}. \label{SMPPROJ}
\end{align}
In practical computation, we adopt the different formulas of the projection $P$
according to $m \le n/2$ or $ m > n/2$. Thus, we give the computational formula
of the projected gradient $Pg_{k}$ as follows:
\begin{align}
       Pg_{k} = \begin{cases}
                     g_{k} - Q_{1}\left(Q_{1}^{T}g_{k}\right), \;
                     \text{if} \; m \le \frac{1}{2}n, \\
                     Q_{2}\left(Q_{2}^{T}g_{k}\right), \; \text{otherwise}.
                \end{cases}
                \label{PROGK}
\end{align}

\vskip 2mm

For a real-world optimization problem \eqref{LEQOPT}, we probably meet the
infeasible initial point $x_{0}$. That is to say, the initial point can not
satisfy the constraint $Ax = b$. We handle this problem by solving the following
projection problem:
\begin{align}
     \min_{x \in \Re^{n}} \; \left\|x - x_{0} \right\|^{2}
     \; \text{subject to} \hskip 2mm  Q_{1}^{T} x = b_{r},  \label{MINDISTVB}
\end{align}
where $b_{r} = \left(R_{1}R_{1}^{T}\right)^{-1}\left(R_{1}b\right)$.
By using the Lagrangian multiplier method and the QR factorization \eqref{ATQR}
of matrix $A^{T}$ to solve problem \eqref{MINDISTVB}, we obtain the initial
feasible point $x_{0}^{F}$ of problem \eqref{LEQOPT} as follows:
\begin{align}
    x_{0}^{F} = x_{0} - Q_{1}\left(Q_{1}^{T}Q_{1}\right)^{-1}
               \left(Q_{1}^{T}x_{0}-b_{r}\right)
              = x_{0} - Q_{1}\left(Q_{1}^{T}x_{0}-b_{r}\right).
    \label{INFIPT}
\end{align}
For convenience, we set $x_{0} = x_{0}^{F}$ in line 4, Algorithm \ref{ALGEPTCTR}.

\vskip 2mm

\subsection{The trusty time-stepping scheme}

\vskip 2mm

Another issue is how to adaptively adjust the time-stepping size $\Delta t_k$
at every iteration. We borrow the adjustment method of the trust-region radius
from the trust-region method due to its robust convergence and fast local
convergence \cite{CGT2000}. After the preprocess of the initial point $x_{0}$,
it is feasible. According to the structure-preserving property of the explicit
continuation method \eqref{EXPTC}-\eqref{LBFGS}, $x_{k+1}$ will preserve the
feasibility. That is to say, $x_{k+1}$ satisfies $Ax_{k+1} = b$. Therefore, we
use the objective function $f(x)$ instead of the nonsmooth penalty function
$f(x) + \sigma \|Ax-b\|$ as the cost function.

\vskip 2mm

When we use the trust-region updating strategy to adaptively adjust time-stepping
size $\Delta t_{k}$ \cite{Higham1999}, we need to construct a local approximation
model of the objective $f(x)$ around $x_{k}$. Here, we adopt the following quadratic
function as its approximation model:
\begin{align}
     q_k(x_{k} + s) = f(x_{k}) + s^{T}g_{k} + \frac{1}{2}s^{T}H_{k}^{-1}s.
     \label{QOAM}
\end{align}
In practical computation, we do not store the matrix $H_{k}$. Thus, we use the
explicit continuation method \eqref{EXPTC}-\eqref{LBFGS} and regard $(H_{k}P)(H_{k}P)^{+}
\approx I$ to simplify the quadratic model $q_k(x_{k}+s_{k}) - q(x_{k})$ as follows:
\begin{align}
    m_{k}(s_{k}) = g_{k}^{T}s_{k}
    - \frac{0.5\Delta t_{k}}{1+\Delta t_{k}}g_{k}^{T}s_{k}
   = \frac{1+0.5\Delta t_{k}}{1+\Delta t_{k}} g_{k}^{T}s_{k}
   \approx q_{k}(x_{k}+s_{k}) - q_{k}(x_{k}).  \label{LOAM}
\end{align}
where $g_{k} = \nabla f(x_k)$. We enlarge or reduce the time-stepping size
$\Delta t_k$ at every iteration according to the following ratio:
\begin{align}
    \rho_k = \frac{f(x_k)-f(x_{k+1})}{m_k(0)-m_k(s_{k})}.
    \label{MRHOK}
\end{align}
A particular adjustment strategy is given as follows:
\begin{align}
     \Delta t_{k+1} = \begin{cases}
          \gamma_1 \Delta t_k, &{if \hskip 1mm 0 \leq \left|1- \rho_k \right| \le \eta_1,}\\
          \Delta t_k, &{if \hskip 1mm \eta_1 < \left|1 - \rho_k \right| < \eta_2,}\\
          \gamma_2 \Delta t_k, &{if \hskip 1mm \left|1-\rho_k \right| \geq \eta_2,}
                   \end{cases} \label{ADTK1}
\end{align}
where the constants are selected as $\eta_1 = 0.25, \; \gamma_1 = 2, \; \eta_2 = 0.75, \;
\gamma_2 = 0.5$  according to numerical experiments. When $\rho_{k} \ge \eta_{a}$,
we accept the trial step $s_{k}$ and let $x_{k+1} = x_{k} + s_{k}$, where
$\eta_{a}$ is a small positive number such as $\eta_{a} = 1.0\times 10^{-6}$.
Otherwise, we discard it and let $x_{k+1} = x_{k}$.

\vskip 2mm

According to the above discussions, we give the detailed implementation of
the explicit continuation method with the trusty time-stepping scheme for the
linearly equality-constrained optimization problem \eqref{LEQOPT} in Algorithm
\ref{ALGEPTCTR}.

\begin{algorithm}
	\renewcommand{\algorithmicrequire}{\textbf{Input:}}
	\renewcommand{\algorithmicensure}{\textbf{Output:}}
    \newcommand{\algorithmicbreak}{\textbf{break}}
    \newcommand{\BREAK}{\STATE \algorithmicbreak}
	\caption{The explicit continuation method with the trusty time-stepping scheme
     for linearly constrained optimization (the Eptctr method)}
    \label{ALGEPTCTR}	
	\begin{algorithmic}[1]
		\REQUIRE ~~\\
        the objective function $f(x)$, the linear constraint $Ax  = b$,  
        the initial point $x_0$ (optional), the terminated parameter 
        $\epsilon$ (optional).
		\ENSURE ~~\\
        the optimal approximation solution $x^{\ast}$.

        \vskip 2mm
        		
        \STATE Set $x_0 = \text{ones}(n, \, 1)$ and $\epsilon = 10^{-6}$ as the default values. 
        \STATE Initialize the parameters: $\eta_{a} = 10^{-6}, \; \eta_1 = 0.25,
        \; \gamma_1 =2, \; \eta_2 = 0.75, \; \gamma_2 = 0.5, \; \theta = 10^{-6}$.
        \STATE Factorize matrix $A^{T}$ with the QR factorization \eqref{ATQR}.
        \STATE Compute
        $$ x_{0} \leftarrow x_{0} - Q_{1}\left(Q_{1}^{T}x_{0}-b_{r}\right),
        $$
        such that $x_{0}$ satisfies the linear system of constraints $Ax = b$.
        \STATE Set $k = 0$. Evaluate $f_0 = f(x_0)$ and $g_0 = \nabla f(x_0)$.
        \STATE Compute the projected gradient $p_{g_{0}} = Pg_{0}$ according to
        the formula \eqref{PROGK}. Set $y_{-1} = 0$ and $s_{-1} = 0$. 
        \STATE Set $\Delta t_0 = 10^{-2}$. 
        \WHILE{$\|p_{g_k}\|> \epsilon$}
          \IF{$\left(|s_{k-1}^{T}y_{k-1}| > \theta s_{k-1}^{T}s_{k-1}\right)$}
            \STATE $ s_{k} = - \frac{\Delta t_{k}}{1 + \Delta t_{k}}
                     \left(p_{g_{k}} - \frac{y_{k-1}(s_{k-1}^{T}p_{g_{k}})
                     + s_{k-1}(y_{k-1}^{T}p_{g_{k}})}{y_{k-1}^{T}s_{k-1}}
                     + 2 \frac{\|y_{k-1}\|^{2}(s_{k-1}^{T}p_{g_{k}})}
                      {(y_{k-1}^{T}s_{k-1})^{2}}s_{k-1}\right)$. 
          \ELSE
            \STATE $s_{k} = - \frac{\Delta t_{k}}{1+ \Delta t_{k}}p_{g_{k}}$.
          \ENDIF
          \STATE Compute $x_{k+1} = x_{k} + s_{k}$.
          \STATE Evaluate $f_{k+1} = f(x_{k+1})$ and compute the ratio $\rho_{k}$
          from equations \eqref{LOAM}-\eqref{MRHOK}.
          \IF{$\rho_k\le \eta_{a}$}
            \STATE Set $x_{k+1} = x_{k}, \; f_{k+1} = f_{k},
            \; p_{g_{k+1}} = p_{g_{k}}, \; g_{k+1} = g_{k}, \; y_{k} = y_{k-1}.$
          \ELSE
            \STATE Evaluate $g_{k+1} = \nabla f(x_{k+1})$.
            \STATE Compute $p_{g_{k+1}} = Pg_{k+1}$ according to the formula
            \eqref{PROGK}. Set $y_{k} = p_{g_{k+1}} - p_{g_{k}}$ and
            $s_{k} = x_{k+1} - x_{k}$.
          \ENDIF
          \STATE Adjust the time-stepping size $\Delta t_{k+1}$ based on the
          trust-region updating scheme \eqref{ADTK1}.
          \STATE Set $k \leftarrow k+1$.
        \ENDWHILE
	\end{algorithmic}
\end{algorithm}

\section{Algorithm Analysis}

In this section, we analyze the global convergence of the explicit continuation
method \eqref{EXPTC}-\eqref{XK1} with the trusty time-stepping scheme and the
L-BFGS updating formula \eqref{LBFGS} for the linearly equality-constrained
optimization problem (i.e. Algorithm \ref{ALGEPTCTR}). Firstly, we give a 
lower-bounded estimate of $m_{k}(0) - m_{k}(s_{k})$ $(k = 1, \, 2, \, \ldots)$. 
This result is similar to that of the trust-region method for the unconstrained
optimization problem \cite{Powell1975}. For simplicity, we assume that the rank of
matrix $A$ is full.

\vskip 2mm

\begin{lemma} \label{LBSOAM}
Assume that the quadratic model $q_{k}(x)$ is defined by equation \eqref{LOAM}
and $s_{k}$ is computed by the explicit continuation method \eqref{EXPTC}-\eqref{LBFGS}.
Then, we have
\begin{align}
    m_{k}(0) - m_{k}(s_{k}) \ge \frac{\Delta t_{k}}{4(1+\Delta t_{k})}
    \left\|p_{g_{k}} \right\|^{2},     \label{PLBREDST}
\end{align}
where $p_{g_k} = Pg_{k} = P\nabla f(x_{k})$ and the projection matrix $P$ is defined
by equation \eqref{PROMAT}.
\end{lemma}
\proof From equation \eqref{LBFGS} and Lemma \ref{LEMHLB}, we know that $H_{k}$
is symmetric positive definite and its eigenvalues are greater than 1/2.
According to the eigenvalue decomposition of $H_{k}$, we know that it exists an
orthogonal matrix $Q_{k}$ such that $H_{k} = Q_{k}^{T}S_{k}Q_{k}$ holds, where
$S_{k} = \text{diag}(\mu_{1}^{k}, \, \ldots, \, \mu_{n}^{k})$ and
$\mu_{i}^{k} \, (i = 1:n)$ are the eigenvalues of $H_{k}$. We denote the
smallest eigenvalue of $H_{k}$ is $\mu_{min}^{k}$. From the explicit
continuation method \eqref{EXPTC} and $P^{2} = P$, we know that
$s_{k} = Ps_{k}$. By combining it with the explicit continuation
method \eqref{EXPTC} and the quadratic model \eqref{LOAM}, we have
\begin{align}
      & m_{k}(0) - m_{k}(s_{k}) \ge -\frac{1}{2}g_{k}^{T}s_{k}
       = -\frac{1}{2}g_{k}^{T}Ps_{k} = - \frac{1}{2}(Pg_{k})^{T}s_{k}
      = \frac{\Delta t_{k}}{2(1+\Delta t_{k})}(Pg_{k})^{T}H_{k}(Pg_{k}) \nonumber \\
      & \quad = \frac{\Delta t_{k}}{2(1+\Delta t_{k})}(Pg_{k})^{T}(Q_{k}^{T}S_{k}Q_{k})(Pg_{k})
      = \frac{\Delta t_{k}}{2(1+\Delta t_{k})}(QPg_{k})^{T}S_{k}(QPg_{k}) \nonumber \\
      & \quad \ge \mu_{min}^{k} \frac{\Delta t_{k}}{2(1+\Delta t_{k})}\|QPg_{k}\|^{2}
      \ge \frac{\Delta t_{k}}{4(1+\Delta t_{k})} \|QPg_{k}\|^{2}
      = \frac{\Delta t_{k}}{4(1+\Delta t_{k})}\|Pg_{k}\|^{2}.  \label{LMOLB}
\end{align}
In the first inequality in equation \eqref{LMOLB}, we use the property
$(1+0.5 \Delta t_{k})/(1+\Delta t_{k}) \ge 0.5$ when $\Delta t_{k} \ge 0$.
Consequently, we prove the result \eqref{PLBREDST}. \qed

\vskip 2mm

In order to prove that $p_{g_k}$ converges to zero when $k$ tends to infinity,
we need to estimate the lower bound of time-stepping sizes
$\Delta t_{k} \, (k = 1, \, 2, \, \ldots)$. We denote the constrained level set
$S_{f}$ as
\begin{align}
     S_{f} = \{x: \; f(x) \le f(x_0), \;  Ax = b \}. \label{LSCONFBD}
\end{align}

\vskip 2mm

\begin{lemma} \label{DTBOUND}
Assume that $f: \; \Re^{n} \rightarrow \Re$ is continuously differentiable and
its gradient $g(x)$ satisfies the following Lipschitz continuity:
\begin{align}
     \|g(x) - g(y)\| \le L_{c} \|x - y\|, \; \forall x, \, y  \in S_{f}.
     \label{LIPSCHCON}
\end{align}
where $L_{c}$ is the Lipschitz constant. We suppose that the sequence $\{x_{k}\}$ is
generated by Algorithm \ref{ALGEPTCTR}. Then, there exists a positive constant
$\delta_{\Delta t}$ such that
\begin{align}
    \Delta t_{k} \ge \gamma_{2} \delta_{\Delta t}     \label{DTGEPN}
\end{align}
holds for all $k = 1, \,  2, \, \dots$, where $\Delta t_{k}$ is adaptively adjusted 
by the trust-region updating scheme \eqref{LOAM}-\eqref{ADTK1}.
\end{lemma}

\vskip 2mm

\proof From Lemma \ref{LEMHLB}, we know that the eigenvalues of $H_{k}$ is greater
than 1/2 and it has at least $n-2$ eigenvalues which equal 1. When $|s_{k-1}^{T}y_{k-1}|
\ge \theta \|s_{k-1}\|^{2}$, we denote the other two eigenvalues of $H_{k}$ as
$\mu_{1}^{k}$ and $\mu_{2}^{k}$. By substituting it into equation
\eqref{ROOTQEQ}, we obtain
\begin{align}
      \mu_{1}^{k} \mu_{2}^{k} = \frac{\|y_{k-1}\|^{2}
      \|s_{k-1}\|^{2}}{(s_{k-1}^{T}y_{k-1})^{2}} \le
      \frac{\|y_{k-1}\|^{2}
      \|s_{k-1}\|^{2}}{\theta^{2} \|s_{k-1}\|^{4}}
      = \frac{1}{\theta^{2}}\frac{\|y_{k-1}\|^{2}}{\|s_{k-1}\|^{2}}.
      \label{TWOEIGS}
\end{align}

\vskip 2mm

From Lemma \ref{LBSOAM} and Algorithm \ref{ALGEPTCTR}, we know $f(x_{k}) \le
f(x_{0}) \, (k = 1, \, 2, \, \ldots)$. From the explicit continuation method
\eqref{EXPTC}-\eqref{LBFGS} and Remark \ref{REMPLS}, we know that $Ax_{k} = Ax_{0} = b
\, (k = 1, \, 2, \, \ldots)$. Thus, from the Lipschitz continuity \eqref{LIPSCHCON}
of $g(x)$, we have
\begin{align}
     \|y_{k-1}\| \le \|P\| \|g(x_{k}) - g(x_{k-1})\| \le L_{C}\|x_{k} - x_{k-1}\|
     = L_{C} \|s_{k-1}\|. \label{UPBYK}
\end{align}
By substituting it into equation \eqref{TWOEIGS} and using $\mu_{i}^{k} > \frac{1}{2}
\, (i = 1, \, 2)$, we obtain
\begin{align}
     \frac{1}{2} < \mu_{i}^{k} < \frac{2L_{C}^{2}}{\theta^{2}}, \; i = 1, \, 2.
     \label{BOUNEIG}
\end{align}
That is to say, the eigenvalues of $H_{k}$ are less than or equal to $M_{H}$,
where $M_{H} = \max \{1, \, {2L_{C}^{2}}/{\theta^{2}}\}$. According to the
eigenvalue decomposition theorem, we know that there exists an orthogonal
matrix $Q_{k}$ such that $H_{k} = Q_{k}^{T}S_{k}Q_{k}$ holds, where
$S_{k} = \text{diag}(\mu_{1}^{k}, \, \ldots, \, \mu_{n}^{k})$ and
$\mu_{i}^{k} \, (i = 1:n)$ are the eigenvalues of $H_{k}$. Consequently, we have
\begin{align}
      \|H_{k}(Pg_{k})\| = \|(Q_{k}^{T}S_{k}Q_{k})Pg_{k}\| =
      \|S_{k}(Q_{k}Pg_{k})\| \le M_{H}\|Pg_{k}\|.  \label{UBHPGK}
\end{align}

\vskip 2mm

From the first-order Taylor expansion, we have
\begin{align}
    f(x_{k}+ s_{k}) =  f(x_{k}) +
    \int_{0}^{1} s_{k}^{T}g (x_{k} + t s_{k}) dt.\label{FOTEFK}
\end{align}
Thus, from equations \eqref{LOAM}-\eqref{PLBREDST}, \eqref{FOTEFK} and the Lipschitz
continuity \eqref{LIPSCHCON} of $g(x)$, we have
\begin{align}
      & \left|\rho_{k} - 1\right| =  \left|\frac{(f(x_{k}) - f(x_{k}+s_{k}))
       - (m_{k}(0) - m_{k}(s_{k}))}{m_{k}(0) - m_{k}(s_{k})}\right|
       \nonumber \\
      & \quad \le \frac{1+\Delta t_{k}}{1+0.5\Delta t_{k}}
      \frac{\left|\int_{0}^{1}s_{k}^{T}(g(x_{k} + t s_{k}) - g(x_{k}))dt\right|}
      {\left|s_{k}^{T}g_{k}\right|} + \frac{0.5\Delta t_{k}}{1+0.5\Delta t_{k}} \nonumber \\
       & \quad \le  \frac{2L_{C}(1+ \Delta t_{k})}{\Delta t_{k}}
       \frac{\|s_{k}\|^{2}}{\|Pg_{k}\|^{2}} + \frac{0.5\Delta t_{k}}{1+0.5\Delta t_{k}}.
       \label{ESTRHOK}
\end{align}
By substituting equation \eqref{EXPTC} and equation \eqref{UBHPGK} into equation
\eqref{ESTRHOK}, we have
\begin{align}
      & \left|\rho_{k} - 1\right|
      \le \frac{2L_{C}\Delta t_{k}}{1+\Delta t_{k}}
       \frac{\|PH_{k}(Pg_{k})\|^{2}}{\|Pg_{k}\|^{2}}
       + \frac{0.5\Delta t_{k}}{1+0.5\Delta t_{k}}       \nonumber \\
      & \le \frac{2L_{C}\Delta t_{k}}{1+\Delta t_{k}}
       \frac{\|P\|^{2} \|H_{k}(Pg_{k})\|^{2}}{\|Pg_{k}\|^{2}}
       + \frac{0.5\Delta t_{k}}{1+0.5\Delta t_{k}}
      \le \frac{(2L_{C} M_{H}^{2}+0.5)\Delta t_{k}}{1+0.5\Delta t_{k}}.
      \label{UBROHK}
\end{align}
In the last inequality of equation \eqref{UBROHK}, we use the property
$\|P\| = 1$. We denote
\begin{align}
     \delta_{\Delta t} \triangleq \frac{\eta_{1}}{2L_{C} M_{H}^{2}+0.5}. \label{UPBMPD}
\end{align}
Then, from equation \eqref{UBROHK}-\eqref{UPBMPD}, when
$\Delta t_{k} \le \delta_{\Delta t}$, it is not difficult to verify
\begin{align}
    \left|\rho_{k} - 1\right| \le (2L_{C} M_{H}^{2}+0.5)\Delta t_{k} \le \eta_{1}.     \label{RHOLETA1}
\end{align}

\vskip 2mm

We assume that $K$ is the first index such that $\Delta t_{K} \le
\delta_{\Delta t}$ where $\delta_{\Delta t}$ is defined by equation \eqref{UPBMPD}.
Then, from equations \eqref{UPBMPD}-\eqref{RHOLETA1}, we know that
$|\rho_{K} - 1 | \le \eta_{1}$. According to the time-stepping adjustment
formula \eqref{ADTK1}, $x_{K} + s_{K}$ will be accepted and the time-stepping size
$\Delta t_{K+1}$ will be enlarged. Consequently, the time-stepping size $\Delta t_{k}$ holds
$\Delta t_{k}\ge \gamma_{2}\delta_{\Delta t}$ for all $k = 1, \, 2, \ldots$. \qed

\vskip 2mm

By using the results of Lemma \ref{LBSOAM} and Lemma \ref{DTBOUND}, we prove
the global convergence of Algorithm \ref{ALGEPTCTR} for the linearly
constrained optimization problem \eqref{LEQOPT} as follows.

\vskip 2mm

\begin{theorem}
Assume that $f: \; \Re^{n} \rightarrow \Re$ is continuously differentiable and
its gradient $\nabla f(x)$ satisfies the Lipschitz continuity \eqref{LIPSCHCON}.
Furthermore, we suppose that $f(x)$ is lower bounded when $x \in S_{f}$, where
the constrained level set $S_{f}$ is defined by equation \eqref{LSCONFBD}.
The sequence $\{x_{k}\}$ is generated by Algorithm \ref{ALGEPTCTR}. Then, we have
\begin{align}
  \lim_{k \to \infty} \inf \|Pg_{k}\| = 0, \label{LIMPGKZ}
\end{align}
where $g_{k} = \nabla f(x_{k})$ and matrix $P$ is defined by equation
\eqref{PROMAT}.
\end{theorem}
\proof According to Lemma  \ref{DTBOUND} and Algorithm \ref{ALGEPTCTR}, we know 
that there exists an infinite subsequence $\{x_{k_{i}}\}$ such that trial steps $s_{k_i}$
are accepted, i.e., $\rho_{k_{i}} \ge \eta_{a}, \, i=1,\, 2,\ldots$. Otherwise,
all steps are rejected after a given iteration index, then the time-stepping
size will keep decreasing, which contradicts \eqref{DTGEPN}. Therefore,
from equations \eqref{MRHOK} and \eqref{PLBREDST}, we have
\begin{align}
     & f(x_{0}) - \lim_{k \to \infty} f(x_{k})
     = \sum_{k = 0}^{\infty} (f(x_{k}) - f(x_{k+1}))
      \nonumber \\
     & \ge \eta_{a} \sum_{i = 0}^{\infty}
    \left(m_{k_{i}}(0) - m_{k_{i}}(s_{k_{i}})\right)
    \ge \eta_{a} \sum_{i = 0}^{\infty}
    \frac{ \Delta t_{k_{i}}}{4(\Delta t_{k_{i}}+1)} \|Pg_{k_{i}}\|.
    \label{LIMSUMFK}
\end{align}

\vskip 2mm

From the result \eqref{DTGEPN} of Lemma \ref{DTBOUND}, we know that
$\Delta t_{k} \ge \gamma_{2} \delta_{\Delta t} \, ( k = 1, \, 2, \, \dots)$.
By substituting it into equation \eqref{LIMSUMFK}, we have
\begin{align}
      f(x_{0}) - \lim_{k \to \infty} f(x_{k}) \ge \eta_{a} \sum_{i = 0}^{\infty}
    \frac{\gamma_{2} \delta_{\Delta t}}
    {4(\gamma_{2} \delta_{\Delta t}+1)}\|Pg_{k_{i}}\|. \label{OBJDMD}
\end{align}
Since $f(x)$ is lower bounded when $x \in S_{f}$ and the sequence $\{f(x_{k})\}$ 
is monotonically decreasing, we have $\lim_{k \to \infty} f(x_{k}) = f^{\ast}$. 
By substituting it into equation \eqref{OBJDMD}, we obtain the result 
\eqref{LIMPGKZ}. \qed

\section{Numerical Experiments}

\vskip 2mm

In this section, some numerical experiments are executed to test the performance
of Algorithm \ref{ALGEPTCTR} (the Eptctr method). The codes are executed by a Dell
G3 notebook with the Intel quad-core CPU and 20Gb memory. We compare Eptctr with
SQP (the built-in subroutine fmincon.m of the MATLAB2018a environment)
\cite{FP1963,Goldfarb1970,MATLAB,NW1999,Wilson1963} Ptctr \cite{LLS2020}) for
some large-scale linearly constrained-equality optimization problems which are
listed in Appendix A. SQP is the traditional and representative optimization
for the constrained optimization problem. Ptctr is significantly better than
SQP for linearly constrained optimization problems according to the numerical
results in \cite{LLS2020}. Therefore, we select these two typical methods as
the basis for comparison.

\vskip 2mm

The termination conditions of the three compared methods are all set by
\begin{align}
    & \|\nabla_{x} L(x_{k}, \, \lambda_{k})\|_{\infty} \le 1.0 \times 10^{-6},
    \label{FOOPTTOL} \\
   & \|Ax_k - b \|_{\infty} \le 1.0 \times 10^{-6}, \;
    k = 1, \, 2, \, \ldots,   \label{FEATOL}
\end{align}
where the Lagrange function $L(x, \, \lambda)$ is defined by equation \eqref{LAGFUN}
and $\lambda$ is defined by equation \eqref{LAMBDA}.

\vskip 2mm

\begin{figure}[!htbp]
      \centering
        \includegraphics[width= 1 \textwidth,height=0.5 \textheight]{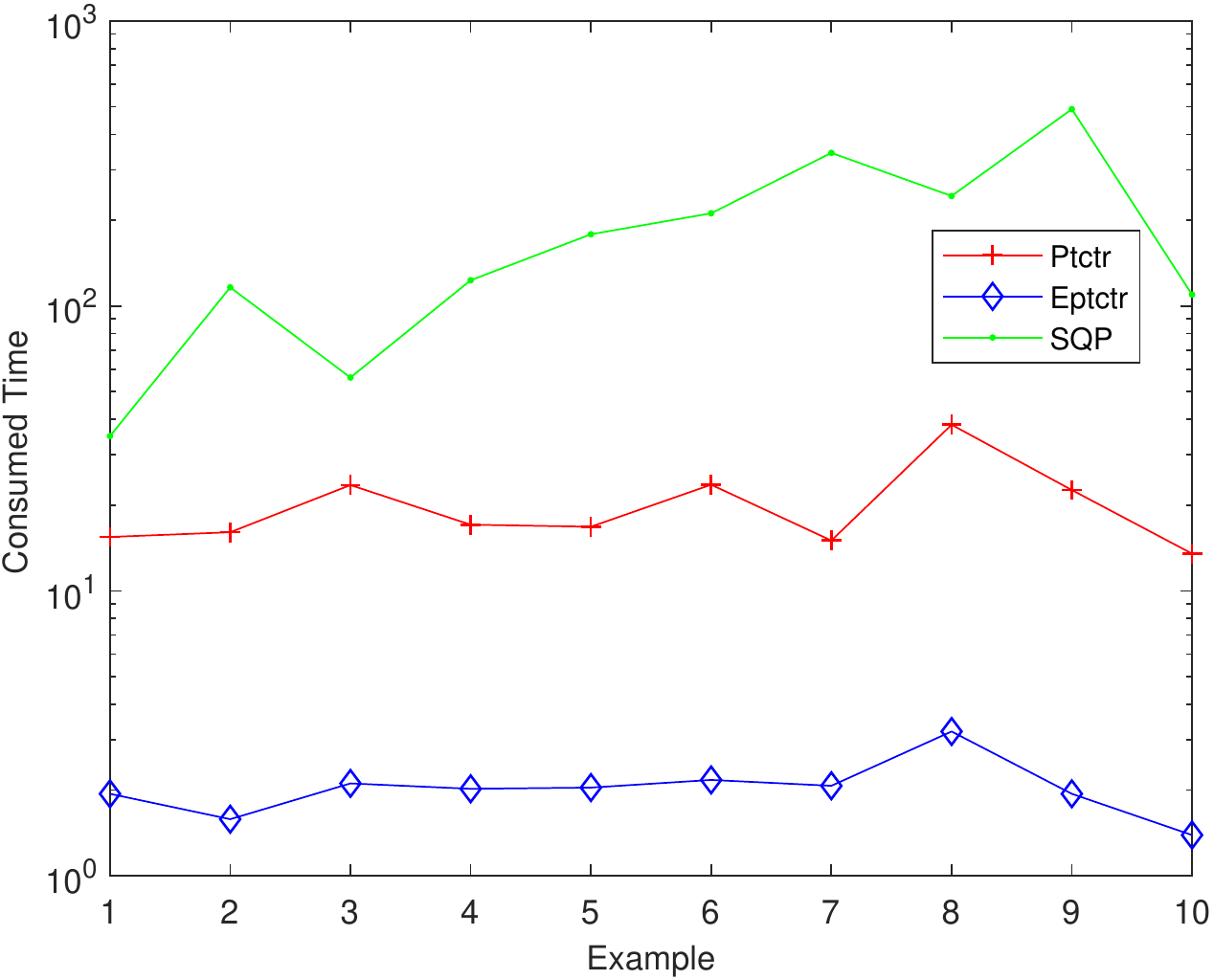}
        \caption{The consumed CPU time (s) of Ptctr, \, Eptctr and SQP for test problems with $n \approx 5000$.}
        \label{fig:CNMTIM}
\end{figure}

We test those ten problems with $n \approx 5000$. The numerical results are 
arranged in Table \ref{TABCOM} and illustrated in Figure \ref{fig:CNMTIM}. 
From Table \ref{TABCOM}, we find that three methods can correctly solve those 
ten test problems and the consumed time of Eptctr is significantly less than 
those of the other two methods for every test problem, respectively. The 
consumed time of Eptctr is about one tenth of that Ptctr or one fifteenth 
to 0.4 percent of that of SQP for the test problem.

\vskip 2mm

From those test data, we find that Eptctr works significantly better
than the other two methods, respectively. One of the reasons is that Eptctr only
involves three pairs of the inner product of two vectors and one matrix-vector product
($p_{g_{k}} = Pg_{k}$) to obtain the trial step $s_{k}$ and involves about $(n-m)n$
flops at every iteration. However, Ptctr needs to solve a linear system of
equations with an $n \times n$ symmetric definite coefficient matrix and involves
about $\frac{1}{3}n^{3}$ flops (p. 169, \cite{GV2013}) at every iteration. SQP
needs to solve a linear system of equations with dimension $(m+n)$ when it solves
a quadratic programming subproblem at every iteration (pp. 531-532, \cite{NW1999})
and involves about $\frac{2}{3}(m+n)^{3}$ flops (p. 116, \cite{GV2013}). Furthermore,
Eptctr can save the storage space of an $(n+m)\times (n+m)$ large-scale matrix, in
comparison to SQP. The required memory of Eptctr is about one fifth of that of SQP.

\vskip 2mm

\begin{table}[!http]
  \newcommand{\tabincell}[2]{\begin{tabular}{@{}#1@{}}#2\end{tabular}}
  \scriptsize
  \centering
  \caption{Numerical results of test problems with $n \approx 5000$.} \label{TABCOM}
  \begin{tabular}{|c|c|c|c|c|c|c|c|c|c|}
  \hline
  \multirow{2}{*}{Problems } & \multicolumn{3}{c|}{Ptctr} & \multicolumn{3}{c|}{Eptctr} & \multicolumn{3}{c|}{SQP}  \\ \cline{2-10}
                        & \tabincell{c}{steps \\(time)}     & $f(x^\star)$  & \tabincell{c}{Mem \\ /Gb} & \tabincell{c}{steps \\(time)}     & $f(x^\star)$  & \tabincell{c}{Mem \\ /Gb}  & \tabincell{c}{steps \\(time)}  & $f(x^\star)$  & \tabincell{c}{Mem \\ /Gb}   \\ \hline

  \tabincell{c}{Exam. 1 \\ (n = 5000, \\ m = n/2)}   & \tabincell{c}{11 \\ (15.46)} &3.64E+04 &3.41   & \tabincell{c}{11 \\ (1.94)} &3.64E+04 &0.51   & \tabincell{c}{2 \\ (34.93)} &3.64E+04 &2.43   \\ \hline

  \tabincell{c}{Exam. 2 \\ (n = 4800, \\ m = n/2)}   & \tabincell{c}{17 \\ (16.06)} &5.78E+03 &4.09   & \tabincell{c}{15 \\ (1.58)} &5.78E+03 &0.31   & \tabincell{c}{14 \\ (116.16)} &5.78E+03 &1.53   \\ \hline

  \tabincell{c}{Exam. 3 \\ (n = 4800, \\ m = 2/3n)}   & \tabincell{c}{12 \\ (23.51)} &2.86E+03 &3.40   & \tabincell{c}{12 \\ (2.11)} &2.86E+03 &0.54   & \tabincell{c}{3 \\ (56.06)} &2.86E+03 &3.08  \\ \hline

  \tabincell{c}{Exam. 4 \\ (n = 5000, \\ m = n/2)}   & \tabincell{c}{11 \\ (17.07)} &493.79 &3.41   & \tabincell{c}{11 \\ (2.02)} &493.79 &0.51   & \tabincell{c}{8 \\ (123.04)} &493.79 &2.43   \\ \hline

  \tabincell{c}{Exam. 5 \\ (n = 5000, \\ m = n/2)}   & \tabincell{c}{14 \\ (16.79)} &432.15 &3.97   & \tabincell{c}{13 \\ (2.04)} &432.15 &0.51   & \tabincell{c}{11 \\ (178.41)} &432.15 &2.43   \\ \hline

  \tabincell{c}{Exam. 6 \\ (n = 4800, \\ m = 2/3n)}   & \tabincell{c}{13 \\ (23.58)} &2.06E+03 &3.57   & \tabincell{c}{15 \\ (2.17)} &2.06E+03 &0.54   & \tabincell{c}{11 \\ (211.35)} &2.06E+03 &3.10   \\ \hline

  \tabincell{c}{Exam. 7 \\ (n = 5000, \\ m = n/2)}   & \tabincell{c}{10 \\ (15.02)} &5.94E+04 &3.22   & \tabincell{c}{13 \\ (2.07)} &5.94E+04 &0.51   & \tabincell{c}{20 \\ (344.37)} &5.94E+04 &2.43   \\ \hline

  \tabincell{c}{Exam. 8 \\ (n = 4800, \\ m = n/3)}   & \tabincell{c}{38 \\ (38.29)} &776.88 &7.36   & \tabincell{c}{133 \\ (3.21)} &-1.21E+04 &0.31   & \tabincell{c}{28 \\ (243.09)} &784.94 &1.53   \\ \hline

  \tabincell{c}{Exam. 9 \\ (n = 5000, \\ m = n/2)}   & \tabincell{c}{12 \\ (22.59)} &2.21E+05 &3.59   & \tabincell{c}{8 \\ (1.94)} &2.21E+05 &0.51   & \tabincell{c}{29 \\ (490.36)} &2.21E+05 &2.43   \\ \hline

  \tabincell{c}{Exam. 10 \\ (n = 4800, \\ m = n/3)}   & \tabincell{c}{16 \\ (13.51)} &2.00 &3.92   & \tabincell{c}{20 \\ (1.39)} &2.00 &0.31   & \tabincell{c}{14 \\ (109.59)} &2.00 &1.53   \\ \hline

\end{tabular}
\end{table}

\vskip 2mm

\vskip 2mm

\vskip 2mm

\section{Conclusion and Future Work}

\vskip 2mm

In this paper, we give an explicit continuation method with the trusty time-stepping
scheme and the L-BFGS updating formula (Eptctr) for linearly equality-constrained
optimization problems. This method only involves three pairs of the inner product of
vector and one matrix-vector product ($p_{g_{k}} = Pg_{k}$) at every iteration,
other than the traditional optimization method such as SQP or the latest continuation
method such as Ptctr \cite{LLS2020}, which needs to solve a quadratic programming
subproblem (SQP) or a linear system of equations (Ptctr). Thus, Eptctr involves
about $(n-m)n$ flops, Ptctr involves about $\frac{1}{3}n^{3}$ flops, and SQP
involves about $\frac{2}{3}(m+n)^{3}$ flops at every iteration. This means that
Eptctr can save much more computational time than SQP or Ptctr. Numerical results
also show that the consumed time of EPtctr is about one tenth of that Ptctr or
one fifteenth to 0.4 percent of that of SQP for the test problem with
$n \approx 5000$. Furthermore, Eptctr can save the storage space of an
$(n+m) \times (n+m)$ large-scale matrix, in comparison to SQP. The required
memory of Eptctr is about one fifth of that of SQP. Therefore, Eptctr is worth
exploring further, and we will extend it to the general nonlinear optimization
problem in the future.

\section*{Acknowledgments} This work was supported in part by Grant 61876199 from National
  Natural Science Foundation of China, Grant YBWL2011085 from Huawei Technologies
  Co., Ltd., and Grant YJCB2011003HI from the Innovation Research Program of Huawei
  Technologies Co., Ltd.. The first author is grateful to Prof. Ya-xiang Yuan and 
  Prof. Li-zhi Liao for their suggestions. 

\begin{appendix}

\section{Test Problems}

\noindent \textbf{Example 1.}
\begin{align}
   \quad  m & = n/2 \nonumber \\
   \min_{x \in \Re^{n}} \;  f(x) & = \sum_{k=1}^{n/2} \;\left(x_{2k-1}^{2} + 10x_{2k}^{2}\right), \;
   \text{subject to} \;  x_{2i-1} + x_{2i} = 4, \; i = 1, \, 2, \ldots, \, m. \nonumber
\end{align}
This problem is extended from the problem of \cite{Kim2010}. We assume that the
feasible initial point is $(2, \, 2, \, \ldots, \, 2, \, 2)$.

\vskip 2mm

\noindent \textbf{Example 2.}
\begin{align}
    \quad   m &= n/3 \nonumber \\
    \min_{x \in \Re^{n}} \;  f(x) & = \sum_{k=1}^{n/2} \; \left(\left(x_{2k-1} -2\right)^{2}
    + 2\left(x_{2k} - 1 \right)^{4}\right) - 5, \;
    \text{subject to} \; x_{3i-2} + 4x_{3i-1}+2x_{3i} = 3, \;  i = 1, \, 2, \ldots, \, n/3.
    \nonumber
\end{align}
We assume that the infeasible initial point is $(-0.5, \, 1.5, \, 1, \, 0, \, \ldots, \, 0, \, 0)$.

\vskip 2mm

\noindent \textbf{Example 3.}
\begin{align}
   \quad  m & = (2/3)n \nonumber \\
   \min_{x \in \Re^{n}} \;  f(x) & = \sum_{k=1}^{n}\; x_{k}^{2}, \;
   \text{subject to} \; x_{3i-2} + 2x_{3i-1} + x_{3i} = 1, \;
     2 x_{3i-2} - x_{3i-1} - 3 x_{3i} = 4, \; i = 1, \, 2, \ldots, \, n/3.  \nonumber
\end{align}
This problem is extended from the problem of \cite{Osborne2016}. The infeasible
initial point is $(1, \, 0.5, \, -1, \, \ldots, \, 1, \, 0.5, \, -1)$.

\vskip 2mm

\noindent \textbf{Example 4.}
\begin{align}
   \quad   m &= n/2 \nonumber \\
   \min_{x \in \Re^{n}} \; f(x) & = \sum_{k=1}^{n/2}\;\left(x_{2k-1}^{2} + x_{2k}^{6}\right) - 1,  \;
   \text{subject to} \;  x_{2i-1} + x_{2i} = 1, \; i = 1, \, 2, \, \ldots, \, n/2. \nonumber
\end{align}
This problem is modified from the problem of \cite{MAK2019}. We assume that the
infeasible initial point is $(1, \, 1, \, \ldots, \, 1)$.

\vskip 2mm

\noindent \textbf{Example 5.}
\begin{align}
   \quad  m & = n/2 \nonumber \\
   \min_{x \in \Re^{n}} \;  f(x) & = \sum_{k=1}^{n/2}\;\left(\left(x_{2k-1} -2\right)^{4}
   + 2\left(x_{2k} -1\right)^{6}\right) - 5, \;
    \text{subject to} \; x_{2i-1} + 4x_{2i} = 3, \; i = 1, \, 2, \, \ldots, \, m. \nonumber
\end{align}
We assume that the feasible initial point is $(-1, \, 1,\, -1, \, 1, \, \ldots, \, -1, \, 1)$.

\vskip 2mm

\noindent \textbf{Example 6.}
\begin{align}
   \quad m & = (2/3)n \nonumber \\
   \min_{x \in \Re^{n}} \;  f(x) &=
   \sum_{k=1}^{n/3}\;\left(x_{3k-2}^{2} + x_{3k-1}^{4} + x_{3k}^{6}\right), \nonumber \\
    \text{subject to} \; &  {x_{3i-2}} + 2{x_{3i-1}} + {x_{3i}} = 1, \;
      2{x_{3i-2}} - {x_{3i-1}} - 3{x_{3i}} = 4, \; i = 1, \, 2, \, \ldots, \, m/2. \nonumber
\end{align}
This problem is extended from the problem of \cite{Osborne2016}. We assume that the
infeasible initial point is $(2, \, 0, \, \ldots, \, 0)$.

\vskip 2mm

\noindent \textbf{Example 7.}
\begin{align}
   \quad m &= n/2 \nonumber \\
   \min_{x \in \Re^{n}} \;  f(x) & = \sum_{k=1}^{n/2}\;\left(x_{2k-1}^{4} + 3x_{2k}^{2}\right), \;
   \text{subject to} \;  x_{2i-1} + x_{2i} = 4, \; i = 1, \, 2, \, \ldots, \, n/2.  \nonumber
\end{align}
This problem is extended from the problem of \cite{Carlberg2009}. We assume that the infeasible initial
point is $(2, \, 2, \, 0, \, \ldots, \, 0, \, 0)$.

\vskip 2mm

\noindent \textbf{Example 8.}
\begin{align}
   \quad  m & = n/3 \nonumber \\
   \min_{x \in \Re{^n}} \;  f(x) & = \sum_{k=1}^{n/3}\;\left(x_{3k-2}^{2} + x_{3k-2}^{2} \, x_{3k}^{2}
   + 2x_{3k-2} \, x_{3k-1} + x_{3k-1}^{4} + 8x_{3k-1}\right), \nonumber \\
   \text{subject to} \; & 2x_{3i-2} + 5x_{3i-1}+x_{3i} = 3, \; i = 1, \, 2, \, \ldots, \, m. \nonumber
\end{align}
We assume that the infeasible initial point is $(1.5, \, 0, \, 0, \, \ldots, \, 0)$.

\vskip 2mm

\noindent \textbf{Example 9.}
\begin{align}
   \quad m &= n/2 \nonumber \\
   \min_{x \in \Re^{n}} \;  f(x) & = \sum_{k =1}^{n/2} \; \left(x_{2k-1}^{4} + 10x_{2k}^{6}\right),
   \;
   \text{subject to} \; x_{2i-1} + x_{2i} =4, \; i = 1, \, 2, \, \ldots, \, m.  \nonumber
\end{align}
This problem is extended from the problem of \cite{Kim2010}. We assume that the feasible
initial point is $(2, \, 2, \, \ldots, \, 2, \, 2)$.

\vskip 2mm

\noindent \textbf{Example 10.}
\begin{align}
   \quad m & = n/3 \nonumber \\
   \min_{x \in \Re^{n}} \;  f(x) & = \sum_{k=1}^{n/3}\;\left(x_{3k-2}^{8} + x_{3k-1}^{6} + x_{3k}^{2}\, \right),
   \;
   \text{subject to} \;  x_{3i-2} + 2x_{3i-1} + 2x_{3i} =1, \; i = 1, \, 2, \, \ldots, \,m. \nonumber
\end{align}
This problem is modified from the problem of \cite{Yamashita1980}.
The feasible initial point is $(1, \, 0, \, 0, \, \ldots, \, 1, \, 0, \, 0)$.

\vskip 2mm

\end{appendix}


\begin{thebibliography}{99}
%
\bibitem{AG2003}
E.~L. Allgower and K. Georg, \emph{Introduction to Numerical
Continuation Methods}, SIAM, Philadelphia, PA, 2003.
%
\bibitem{AP1998}
U.~M. Ascher and L.~R. Petzold, \emph{Computer Methods for Ordinary
Differential Equations and Differential-Algebraic Equations}, SIAM,
Philadelphia, PA, 1998.
%
\bibitem{Bertsekas2018}
D.~P. Bertsekas, \emph{Nonlinear Programming (3rd Edition)}, Tsinghua University
Press, 2018.
%
\bibitem{BB1989}
A.~A. Brown and M.~C. Bartholomew-Biggs, \emph{ODE versus SQP
methods for constrained optimization}, Journal of Optimization and
Theory Applications, \textbf{62} (3): 371-386, 1989.
%
\bibitem{BCP1996}
K.~E. Brenan, S.~L. Campbell and L.~R. Petzold, \emph{Numerical
solution of initial-value problems in differential-algebraic
equations}, SIAM, Philadelphia, PA, 1996.
%
\bibitem{BNY1987}
R. Byrd, J. Nocedal and Y.~X. Yuan, \emph{Global convergence of a class of
quasi-Newton methods on convex problems}, SIAM Journal of Numerical Analysis,
\textbf{24}: 1171-1189, 1987.
%
\bibitem{Carlberg2009}
K. Carlberg, \emph{Lecture notes of constrained optimization},
\url{https://www.sandia.gov/~ktcarlb/opt_class/OPT_Lecture3.pdf}, 2009.
%
\bibitem{CMFO2009}
F. Caballero, L. Merino, J. Ferruz and A. Ollero, \emph{Vision-based odometry
and SLAM for medium and high altitude flying UAVs}, Journal of Intelligent and
Robotic Systems, \textbf{54} (1-3): 137-161, 2009.
%
\bibitem{CKK2003}
T.~S. Coffey, C.~T. Kelley and D.~E. Keyes, \emph{Pseudotransient
continuation and differential-algebraic equations}, SIAM Journal on
Scientific Computing, \textbf{25}: 553-569, 2003.
%
\bibitem{CGT2000}
A.~R. Conn, N. Gould and Ph.~L. Toint, \emph{Trust-Region Methods},
SIAM, Philadelphia, USA, 2000.
%
\bibitem{FM1990}
A.~V. Fiacco and G.~P. McCormick, \emph{Nonlinear programming: Sequential
Unconstrained Minimization Techniques}, SIAM, 1990.
%
\bibitem{FP1963}
R. Fletcher and M.~J.~D. Powell, \emph{A rapidly convergent descent method for
minimization}, Computer Journal, \textbf{6}: 163-168, 1963.
%
\bibitem{Goh2011}
B.~S. Goh, \emph{Approximate greatest descent methods for optimization with
equality constraints}, Journal of Optimization Theory and Applications
\textbf{148} (3): 505-527, 2011.
%
\bibitem{Goldfarb1970}
D. Goldfarb, \emph{A family of variable metric updates derived by variational means},
Mathematics of Computing, \textbf{24}: 23-26, 1970.
%
\bibitem{GV2013}
G.~H. Golub and C.~F. Van Loan, \emph{Matrix Computations}, 4th ed., The Johns
Hopkins University Press, 2013.
%
\bibitem{HLW2006}
E. Hairer, C. Lubich and G. Wanner, \emph{Geometric Numerical Integration:
Structure-Preserving Algorithms for Ordinary Differential Equations}, 2nd ed.,
Springer, Berlin, 2006.
%
\bibitem{HM1996}
U. Helmke and J.~B. Moore, \emph{Optimization and Dynamical Systems}, 2nd ed.,
Springer-Verlag, London, 1996.
%
\bibitem{Higham1999}
D.~J. Higham, \emph{Trust region algorithms and timestep selection}, SIAM Journal
on Numerical Analysis, \textbf{37}: 194-210, 1999.
%
\bibitem{KLQCRW2008}
C.~T. Kelley, L.-Z. Liao, L. Qi, M.~T. Chu, J.~P. Reese and C.
Winton, \emph{Projected Pseudotransient Continuation}, SIAM Journal
on Numerical Analysis, \textbf{46}: 3071-3083, 2008.
%
\bibitem{LF2000}
D.~G. Liu and J.~G. Fei, \emph{Digital Simulation Algorithms for Dynamic Systems}
(in Chinese), Science Press, Beijing, 2000.
%
\bibitem{LL2010}
S.-T. Liu and X.-L. Luo, \emph{A method based on Rayleigh quotient gradient
flow for extreme and interior eigenvalue problems}, Linear Algebra and its
Applications, \textbf{432} (7): 1851-1863, 2010.
%
\bibitem{Luo2012}
X.-L. Luo, \emph{A dynamical method of DAEs for the smallest eigenvalue
problem}, Journal of Computational Science, \textbf{3} (3): 113-119, 2012.
%
\bibitem{LKLT2009}
X.-L. Luo, C.~T. Kelley, L.-Z. Liao and H.-W. Tam, \emph{Combining trust-region
techniques and Rosenbrock methods to compute stationary points},
Journal of Optimization Theory and Applications, \textbf{140} (2): 265-286, 2009.
%
\bibitem{LLW2013}
X.-L. Luo, J.-R. Lin and W.-L. Wu, \emph{A prediction-correction dynamic
method for large-scale generalized eigenvalue problems}, Abstract and
Applied Analysis, Article ID 845459, 1-8,
\url{http://dx.doi.org/10.1155/2013/845459}, 2013.
%
\bibitem{LLS2020}
X.-L. Luo, J.-H. Lv and G. Sun, \emph{Continuation method with the trusty time-stepping
scheme for linearly constrained optimization with noisy data}, published online in 
\url{http://arxiv.org/abs/2005.05965} or \url{http://doi.org/10.1007/s11081-020-09590-z}, 
Optimization and Engineering, Accepted, January 10, 2021.
%
\bibitem{LXL2020}
X.-L. Luo, H. Xiao and J.-H. Lv, \emph{Continuation Newton methods with the
residual trust-region time-stepping scheme for nonlinear equations}, June 2020,
arXiv preprint, \url{http://arxiv.org/abs/2006.02634}.
%
\bibitem{LY2021}
X.-L. Luo and Y.Y. Yao, \emph{Primal-dual path-following methods and the
trust-region updating strategy for linear programming with noisy data}, June 2020, arXiv
preprint available at \url{http://arxiv.org/abs/2006.07568}, minor revision 
resubmitted to Journal of Computational Mathematics, January 16, 2021.
%
\bibitem{MAK2019}
M.-W. Mak, \emph{Lecture notes of constrained optimization and support vector machines},
\url{http://www.eie.polyu.edu.hk/~mwmak/EIE6207/ContOpt-SVM-beamer.pdf}, 2019.
%
\bibitem{MATLAB}
MATLAB 9.4.0 (R2018a), The MathWorks Inc., \url{http://www.mathworks.com}, 2018.
%
\bibitem{NW1999}
J. Nocedal and S.~J. Wright, \emph{Numerical Optimization},
Springer-Verlag, 1999.
%
\bibitem{Kim2010}
N.~H. Kim, \emph{Leture notes of constrained optimization},
\url{https://mae.ufl.edu/nkim/eas6939/ConstrainedOpt.pdf}, 2010.
%
\bibitem{Osborne2016}
M.~J. Obsborne, \emph{Mathematical methods for economic theory},
\url{https://mjo.osborne.economics.utoronto.ca/index.php/tutorial/index/1/mem}, 2016.
%
\bibitem{Pan1992}
P.-Q. Pan, \emph{New ODE methods for equality constrained optimization (2):
algorithms}, Journal of Computational Mathematics, \textbf{10} (2): 129-146, 1992.
%
\bibitem {Powell1975}
M.~J.~D. Powell, \emph{Convergence properties of a class of minimization
algorithms}, in: O.L. Mangasarian, R. R. Meyer and S. M. Robinson,
eds., \emph{Nonlinear Programming 2}, Academic Press, New York, 1-27, 1975.
%
\bibitem{Schropp2000}
J. Schropp, \emph{A dynamical systems approach to constrained minimization},
Numerical Functional Analysis and Optimization, \textbf{21} (3-4): 537-551, 2000.
%
\bibitem{Schropp2003}
J. Schropp, \emph{One- and multistep discretizations of index 2
differential algebraic systems and their use in optimization},
Journal of Computational and Applied Mathematics, \textbf{150}: 375-396, 2003.
%
\bibitem{SGT2003}
L.~F. Shampine, I. Gladwell and S. Thompson, \emph{Solving ODEs with MATLAB},
Cambridge University Press, Cambridge, 2003.
%
\bibitem{SY2006}
W.~Y. Sun and Y.~X. Yuan, \emph{Optimization Theory and Methods: Nonlinear
Programming}, Springer, New York, 2006.
%
\bibitem{Tanabe1980}
K. Tanabe, \emph{A geometric method in nonlinear programming}, Journal of
Optimization Theory and Applications, \textbf{30} (2): 181-210, 1980.
%
\bibitem{Simos2013}
T.~E. Simos, \emph{New open modified Newton Cotes type formulae as multilayer 
symplectic integrators}, Applied Mathematical Modelling \textbf{37}: 1983-1991, 
2013. 
%
\bibitem{USS2020}
N. Ullah, J. Sabi and A. Shah, \emph{A derivative-free scaled memoryless BFGS
method for solving a system of monotone nonlinear equations}, Submitted to Numerical
Linear Algebra with Applications, October 2020.
%
\bibitem{Wilson1963}
R.~B. Wilson, \emph{A Simplicial Method for Convex Programming}, Phd thesis,
Harvard University, 1963.
%
\bibitem{Yamashita1980}
H. Yamashita, \emph{A differential equation approach to nonlinear programming},
Mathematical Programming, \textbf{18}: 155-168, \url{https://doi.org/10.1007/BF01588311},
1980.
%
\bibitem{Yuan2015}
Y. Yuan, \emph{Recent advances in trust region algorithms}, Mathematical Programming,
\textbf{151}: 249-281, 2015.
%
%\bibitem{ZS2015}
%J. Zhang and S. Singh, \emph{Visual-inertial combined odometry system
%for aerial vehicles}, Journal of Field Robotics, \textbf{32} (8): 1043-1055, 2015.

\end{thebibliography}
\end{document}